\documentclass[11pt]{article}

\usepackage{amsmath}
\usepackage{amsthm}
\usepackage{amssymb}
\usepackage[left=2.5cm,right=2.5cm,top=2.5cm,bottom=2.5cm]{geometry}
\usepackage{tikz}
\usetikzlibrary{positioning,shapes,calc,backgrounds,fit,matrix}
\usepackage{pgfplots}
\usetikzlibrary{plotmarks}
\usepackage{multirow}

\usepackage{algorithm}
\usepackage{algcompatible}

\long\def\drop#1{}

\begin{document}

\title{Multilevel Monte Carlo Simulations of Composite Structures with Uncertain Manufacturing Defects}

\author{
T. J. Dodwell\footnote{College of Engineering, Mathematics and
  Physical Sciences, University of Exeter, UK \& The Alan Turing
  Institute, London, NW1 2DB, UK. Email: t.dodwell@exeter.ac.uk}, 
S. Kynaston\footnote{Department of
  Mathematical Sciences, University of Bath, Bath, UK.},
R. Butler\footnote{Department of Mechanical Engineering, University of
  Bath, Bath BA2 7AY, UK.}, R. T. Haftka\footnote{Department of
  Mechanical and Aerospace Engineering, University of Florida, FL,
  USA.}, Nam H. Kim$^{\S}$ and R. Scheichl$^{\dag,}$\footnote{Interdisciplinary
  Center for Scientific Computing, University of Heidelberg, Germany}
}

\maketitle
\noindent

\begin{abstract}
By adopting a Multilevel Monte Carlo (MLMC) framework, we show that only a handful of costly fine scale computations are needed to accurately estimate statistics of the failure of a composite structure, as opposed to the thousands typically needed in classical Monte Carlo analyses. We introduce the MLMC method, compare its theoretical complexity with classical Monte Carlo, and give a simple-to-implement algorithm which includes a simple extension called MLMC with selective refinement to efficiently calculated structural failure probabilities. To demonstrate the huge computational gains we present two benchmark problems in composites: (1) the effects of fibre waviness on the compressive strength of a composite material, (2) uncertain buckling performance of a composite panel with uncertain ply orientations. For our most challenging test case, estimating a rare ($\sim 1/150$) probability of buckling failure of a composite panel, we see a speed-up factor $> 1000$. Our approach distributed over $1024$ processors reduces the computation time from $218$ days to just $4.5$ hours. This level of speed up makes stochastic simulations that would otherwise be unthinkable now possible.
\end{abstract}

\section{Introduction}

Within the aerospace manufacturing sector, where safety is paramount, risk
is quantified and reduced by heuristic safety factors and
expensive programmes of empirical testing over a variety of length scales
before a new designs
can enter production, with more tests at coupon than at component scale,
the so-called test pyramid.
The high cost of certification and the inefficiency of general safety factors
has led to new initiatives~\cite{AC20}
whereby numerical simulation and stochastic methods and an increasing interest in probabilistic design \cite{Aca07}. Both of which provided opportunities to demonstrate
structural integrity even when experimental / statistical data is incomplete,  offering scope to challenge conservative failure
limits and reduce design-to-manufacture time.

\smallskip
In complex composite manufacturing
processes, uncertainty arises from a number of different sources,
e.g. material variability \cite{Sut12}, machine tolerance \cite{Rhe13} and
process-induced defects such as fibre waviness or ply wrinkling
\cite{Dod14,Fle15,Bel18,San18}. However, statistical simulations typically require a
large number of analyses, and thus can become extremely
computationally expensive. For that reason a host of techniques for mitigating their
cost has been developed in the engineering and statistical communities.
When the probability of failure due to a single failure mode is all
that is needed, methods, such as the first order reliability method (FORM)
\cite[Section~4.4]{Mel99}, can substantially reduce the number of
required simulations. Such methods have been applied to buckling of
shells with random imperfections \cite{Eli87} and composite laminates \cite{Sha10}. However, such
methods do not capture the statistical interaction between multiple
failure modes; and degenerate further when random variables have complex non-Gaussian distributions and multiple Most Probable Points \cite{Par15}. In such cases, Monte Carlo simulations  are often the only choice for
capturing such interactions.
There are therefore a host of methods that are targeted at reducing
the cost of Monte Carlo simulations. Importance sampling methods,
e.g. \cite[Section~3.4]{Mel99},  reduce the number of required
simulations by preferential sampling near the boundary that is 
separating the safe and the failure domains. However, to identify this
boundary can be as difficult as the original Monte Carlo
simulation. Similarly, separable Monte Carlo
methods \cite{Sma10} take advantage of the independence of
uncertainty sources, and the two approaches can be combined for
additional savings \cite{Cha13}.
Surrogates are often used for allowing large Monte Carlo sampling 
\cite{allaire}; however, surrogates suffer from the `curse of
dimensionality'. One approach for alleviating this
problem is to combine a large number of low-fidelity, inexpensive
simulations with a small number of higher fidelity simulations. For
example, Alexandrov et al. (2001) describes the use of multiple model resolutions
for constructing  surrogates for aerodynamic optimisation \cite{Ale01}.

\smallskip
In our application, where a large number of defects need to be
simulated it would be impractical to construct accurate surrogates.
However, we can still take advantage of combining fidelities with
different mesh sizes. This paper, therefore, sets out to optimise the
use of a hierarchy of coarse and fine finite element (FE) models for 
Monte-Carlo simulations of composites with defects. By adopting a Multilevel Monte Carlo (MLMC) framework, we show that only 
a handful of costly fine scale computations are needed to accurately 
estimate statistics of structural failure loads, as opposed 
to thousands of fine scale samples typically needed in classical Monte 
Carlo analyses. The missing exploration of the variability, leading to sampling error, is taken 
care of by a large number of coarse simulations. 
Multilevel techniques were first suggested in the 
context of option pricing in financial mathematics \cite{Gil08}. 
Its huge potential in uncertainty quantification for engineering 
applications was identified by Cliffe et al. \cite{Cli11} where 
it has been motivated via a subsurface hydrology application. 
Since then it has been applied to a range of other 
applications \cite{barth,mueller,mishra}, it has been 
improved \cite{collier,Elf16} and extended to allow also for 
experimental data to be taken into account in a Bayesian 
setting \cite{hoang,Dod15}. 

\smallskip
Importantly in many engineering applications, estimating the expected load of structural failure is of limited interest, instead often in design we wish to compute the probability
that the failure load is less than a `safe' value. Such a model has a binary output, failure ($1$) or not ($0$). In this paper we propose an extension to MLMC, motivated from an approach proposed by Elferson et al. \cite{Elf16}. By using an error estimator, it is possible for most samples to conclude from a coarse, computationally cheap, model that further model refinements will also not fail, providing that the coarse model predicts a load sufficiently far from the failure boundary. We refer to this extension as Multilevel Monte Carlo with Selective Refinement (MLMC-SR), and demonstrate that it delivers significant further computational gains over even MLMC.

\smallskip
In this paper we describe the multilevel Monte Carlo method in 
a fairly abstract way (Sec. \ref{sec:MLMC}) to show that it can be applied to a broad class of problems in composite applications. We compare its theoretical complexity with that of a standard Monte Carlo simulation and provide simple-to-implement, practical algorithms for both MLMC and MLMC-SR.
To demonstrate the huge computational gains and theoretical results we present two benchmark/classical analysis problems in aerospace composites. For the first we explore the effects of fibre waviness on the compressive strength of a composite material, and for the second test we consider the buckling performance of a skin panel with uncertain ply orientations. The numerical experiments in Section \ref{sec:results} confirm the 
theoretically predicted gains for the model problems with huge potential speed-ups as much as $1000$ fold. This level of speed-up brings stochastic simulations that would otherwise be unthinkable into the
feasible range.

\smallskip
From an engineering viewpoint, whilst the model problems are chosen to represent the
typical gains that can be achieved with the MLMC methodology, in addition, we also learn
something about the engineering implications of uncertainty in each case. In the buckling test problem, perhaps unsurprisingly, the numerical results show that random
variations in ply angles increase the risk of buckling failure significantly. With variations in ply angles of the order typically observed in an Automated Fibre Placement (AFP) machine ($\pm 5^\circ$) significant variability is observed in buckling performance. As for our numerical study of the effects of random fibre waviness on the compressive strength of composites, high fidelity stochastic simulations show remarkably good agreement with Budiansky's classical kinking model \cite{Bud83} if the misalignment angle is taken to be standard deviation of the misalignment random field.

\section{Multi-Level Monte Carlo Methodology and Implementation }\label{sec:MLMC}

To describe the multilevel uncertainty quantification method, let us assume we have a finite element model of a composite structure that is subject to some uncertainty in its material properties, for example due to a defect or the misalignment of fibres. The accuracy and the computational cost of the model is directly linked to the number of degrees of freedom ($M$) and thus to the resolution of the finite element mesh. Typically, for a particular application we are interested in some scalar quantity of interest $Q$. This may be point values of finite element solution (i.e. displacement), or a more complicated nonlinear functional (e.g. failure stress). In the context of the example problems we consider, it is the expected value of a failure stress or a critical buckling load. In cases with random defects or uncertainty, we are therefore interested in estimating the expected value of $Q$, denoted $\mathbb E[Q]$, or perhaps the distribution of $Q$.

\subsection{Standard Monte-Carlo Simulation}

In a typical Monte Carlo (MC) analysis, we create a large number ($N$) of independent random realisations (or samples) of our parameters. For each sample we compute the FE solution on an mesh with $M$ degrees of freedom. From this solution the quantity of interest of the $j^{th}$ sample, $Q^{(j)}_M$, is computed. The average 
\begin{equation}\label{eq:MC_estimator}
\widehat{Q}^{\text{MC}}_{M,N} = \frac{1}{N} \sum_{j=1}^N Q_M^{(j)}
\end{equation}
of these independent samples of $Q_M$ is then the standard Monte Carlo estimator for the expected value $\mathbb{E}[Q_M]$ of $Q_M$.

\smallskip\noindent
The total error is quantified via the \textit{root mean square error} (RMSE), given by
\begin{equation}
\label{eq:MC_MSE_1}
e(\widehat{Q}^{\text{MC}}_{M,N}) = \left(\mathbb{E}[(\widehat{Q}^{\text{MC}}_{M,N} - \mathbb{E}[Q])^2]\right)^{1/2}.
\end{equation}
The mean square error can be expanded so that
\begin{equation}\label{eq:MC_MSE_2}
e(\widehat{Q}_{M,N}^{\text{MC}})^2 = \mathbb{E}[Q_M
- Q]^2 +  \frac{\mathbb{V}[Q_M]}{N}.
\end{equation}
where $\mathbb{V}[Q_M]$ denotes the variance of the random variable $Q_M$.
From the expression we can identify two sources of error in the estimator \eqref{eq:MC_estimator}. The first term is the square of the {\em bias error}. This arises since we are actually interested in the expected value $\mathbb{E}[Q]$ of $Q$, the (inaccessible) random variable corresponding to the exact solution without any FE error. However, since the FE method converges for each sample, as $M \to \infty$, we also have
\begin{equation}\label{eq:M_alpha}
|\mathbb{E}[Q_M-Q]| \leq  C_1 M^{-\alpha},
\end{equation} 
where $\alpha>0$ is the order of convergence, and $C_1$ is some constant independent of $M$. We can reduce this error below any prescribed bias tolerance $e_b$ by making $M$ sufficiently large.  

\smallskip\noindent
The second term of \eqref{eq:MC_MSE_2} gives the {\em sampling error} since we only approximate $\mathbb{E}[Q]$ with $N$ samples.  To ensure this term  is smaller than a sample tolerance $e_s^2$, it suffices to choose 
\begin{equation}
N \ge \frac{V}{e_s^{2}},
\end{equation}
where $V = \mathbb{V}[Q] \approx \mathbb{V}[Q_M]$, for $M$ sufficiently large. The total mean square error is then less than
$e_b^2 + e_s^2$. To ensure that this is less than $e^2$ we can choose 
\begin{equation}\label{tausplit}
e_s^2 = \theta e^2 \quad \text{and} \quad e_b^2 = (1-\theta)e^2, \quad \text{for some} \ \  0 < \theta < 1.
\end{equation}

\smallskip\noindent
We observe that in order to reduce the total error in \eqref{eq:MC_MSE_1} it is necessary to increase both the number of degrees of freedom $M$ and the number of samples $N$. This very quickly leads to an intractable problem when the cost to compute each sample to a sufficiently high accuracy is high. The cost $\mathcal C$ for one sample $Q_M^{(j)}$ of $Q_M$, in terms of floating
point operations (FLOPs) or CPU time, depends on the complexity of the FE solver. Typically it will grow
like 
\begin{equation}\label{eqn:gamma}
\mathcal C_\ell \leq C_2 M_\ell^\gamma,
\end{equation}
for some $\gamma \ge 1$ and constant $C_2$, independent of both $j$ and $M$. Thus, the total cost to achieve a root mean square error $e(\widehat{Q}_{M,N}^{\text{MC}}) \le e$ with standard MC is
\begin{equation}
\text{Cost}(\widehat{Q}_{M,N}^\text{MC}) \geq C_2 N M^\gamma \geq  \, C_3 e^{-2-\gamma/\alpha}.
\end{equation}

 
\subsection{Multilevel Monte-Carlo Simulation (MLMC)}\label{subsec:MLMC} 
\begin{figure}
\centering
\includegraphics[width = \linewidth]{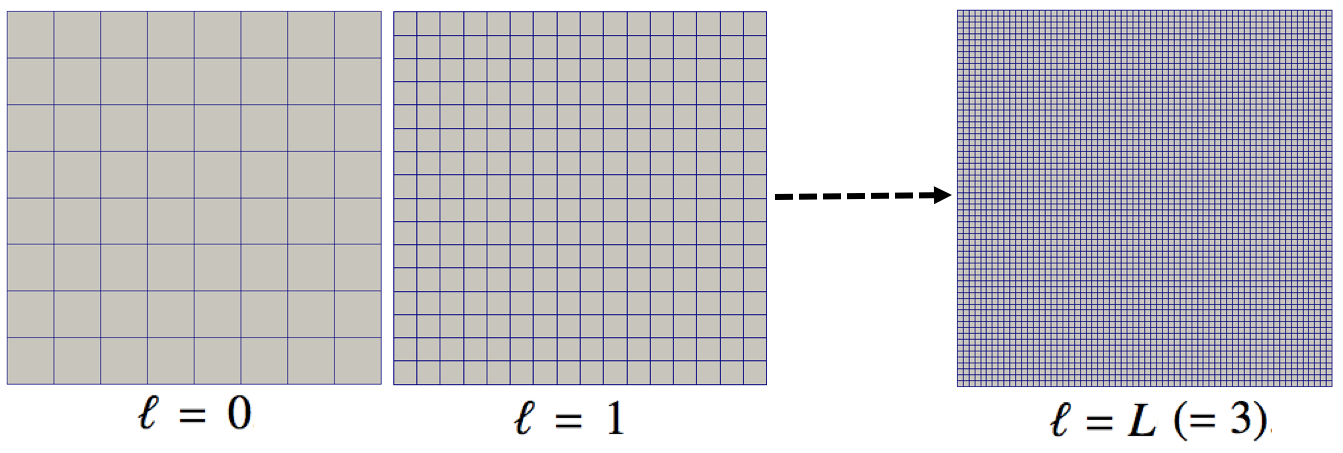}
\caption{Example hierarchy of two-dimensional, quadrilateral finite element meshes for the multilevel algorithm achieved through uniform refinement.}
\label{fig:meshes}
\end{figure}

Multilevel Monte Carlo simulation (MLMC) \cite{Gil08,Cli11} seeks to reduce the variance of the estimator \eqref{eq:MC_estimator} and thus to reduce computational time, by recursively using a hierarchy of FE models as control variants. The standard MC estimator in the previous section was too costly because all samples were computed to the required level to sufficiently reduce the discretisation (or bias) error. Let us now introduce a hierarchy of FE models, obtained by refinement of a coarse mesh as shown in Fig.~\ref{fig:meshes}. Each mesh corresponds to a {\em level} $0 \le \ell \le L$ in our multilevel method with $M_0 < \dots < M_\ell < \dots < M_L$ degrees of freedom, respectively, where $M_0$ is typically small.

\smallskip \noindent
By exploiting the linearity of the expectation operator, the MLMC method avoids estimating $\mathbb{E}[Q]$ directly on the finest, most computationally expensive, level $L$. Instead it estimates the mean on the coarsest level, and corrects
this mean successively by adding estimates of the expected values of differences between subsequent levels, $Y_\ell = Q_{M_\ell} - Q_{M_{\ell-1}}$, for $\ell \ge 1$; i.e. using the identity
\begin{equation}\label{eq:MLMC_expectation}
\mathbb{E}[Q_M] = \mathbb{E}[Q_{M_0}] + \sum_{\ell=1}^L\mathbb{E}[Y_\ell]\,.
\end{equation}
The MLMC estimator for $\mathbb E[Q]$ is then given by
\begin{equation}\label{eq:MLMC_defn}
\widehat{Q}_M^{\text{ML}} = \widehat{Q}^{\text{MC}}_{M_0,N_0} + \sum_{\ell=1}^{L}\widehat{Y}^{\text{MC}}_{\ell,N_\ell}
\end{equation}
where the numbers of samples $N_\ell$ are judiciously chosen to minimise the total cost of this estimator for a given prescribed sampling error (see below). Note that samples $Y^{(j)}_\ell$ of $Y_\ell$ require the FE approximations $Q^{(j)}_{M_\ell}$ and $Q^{(j)}_{M_{\ell-1}}$ on two consecutive mesh levels, i.e. two solves, but crucially both with the same sample $\boldsymbol{\xi}^{(j)}$ of the parameters. 

\medskip\noindent
The cost of the MLMC estimator is
\begin{equation}
\label{MLMC_cost}
\text{Cost}(\widehat{Q}_{M_L}^\text{ML}) =  \sum_{\ell=0}^{L} N_\ell
\mathcal{C}_\ell\,,
\end{equation}
where $\mathcal{C}_\ell$ is the cost to compute a single sample of $Y_\ell$ (resp.~$Q_{M_0}$) on level $\ell$ (resp.~0). By using independent samples across all the levels, the mean square error of $\widehat{Q}_M^{\text{ML}}$ expands to
\begin{equation}
e(\widehat{Q}_M^{\text{ML}})^2 = \big(\mathbb{E}[Q_M-Q]\big)^2 \; + \;
\sum_{l=0}^L N_\ell^{-1}V_\ell\,,
\end{equation}
where $V_0 = \mathbb{V}[Q_{M_0}]$ and $V_\ell = \mathbb{V}[Y_\ell]$, for
$\ell\ge 1$. This leads to a hugely reduced
variance of the estimator since both FE solutions $Q_{M_\ell}$ and $Q_{M_{\ell-1}}$
converge to $Q$ and thus
$$
V_\ell = \mathbb{V}[Q_{M_\ell}-Q_{M_{\ell-1}}]\to 0 \quad \text{as}
\quad M_{\ell} \to
\infty.
$$ 
Let us assume that 
\begin{equation}\label{eqn:beta}
V_\ell \leq C_4 M_\ell^{-\beta}.
\end{equation}
As for the standard MC estimator, we can ensure that the bias error is less than $e_b$ by choosing $M = M_L$ sufficiently large to satisfy \eqref{eq:M_alpha}. To choose the numbers of samples $N_\ell$ on each of the levels and thus to ensure that the sampling error is less than $e_s$, we still have some freedom and we will use this to minimise the computational cost of the overall MLMC algorithm. The samples per level are chosen by a constrained optimisation problem which minimises $\text{Cost}(\widehat{Q}_{M}^\text{ML})$ \eqref{MLMC_cost} with respect to $N_0, \ldots, N_\ell$, subject to the constraint that the samples over all levels are sufficient to reduce the sampling error of the multilevel estimate \eqref{eq:MLMC_expectation} below the required tolerance, such that $$\sum_{\ell=0}^L N_\ell^{-1}V_\ell = e_s^2.$$ This leads to
\begin{equation}
\label{eqn:optimal_Nl}
N_\ell \ = \ e_s^{-2} \left(\sum_{\ell=0}^L \sqrt{V_\ell
    \mathcal{C}_\ell} \right) \, \sqrt{\frac{V_\ell}{C_\ell}}
\end{equation}
For which the total cost using \eqref{MLMC_cost} is
\begin{equation}\label{eqn:mlmcComplexity}
\text{Cost}(\widehat{Q}_{M}^\text{ML}) \; = \; e^{-2} \left(\sum_{\ell=0}^L \sqrt{V_\ell
    \mathcal{C}_\ell} \right)^2 \; \le \; C_5 \, e^{-2-\max\left(0,\frac{\gamma-\beta}{\alpha}\right)},
\end{equation}
where $\alpha,\beta$ and $\gamma$ are as defined above and $e$ is again the tolerance for the total root mean square error.

\smallskip
\noindent
There are three regimes which determine the computational cost of a MLMC algorithm:
\begin{enumerate}
\item \smallskip If the variance $V_\ell$ decays faster than
the cost $\mathcal{C}_\ell$ grows (with respect to $\ell$), i.e $\beta > \gamma$, then the majority of the work is on
level $0$ and the total cost is proportional to $e^{-2}$
\item \smallskip 
If $V_\ell$ decays slower than $\mathcal{C}_\ell$ grows, i.e $\beta < \gamma$ then the majority of the work is on
level $L$ and the total cost is proportional to
$e^{-2-\frac{\gamma-\beta}{\alpha}}$
\item
\smallskip 
If $V_\ell C_\ell$ is bounded, i.e. $\beta = \gamma$, then the work is spread evenly over all levels and $C_5$ grows with $(\log e)^2$. 
\end{enumerate}

\smallskip In the work presented in this paper we consider a hierarchy of levels created by a uniform refinement of a coarse mesh. In this case for each random sample, the mesh at a given level is identical. Recent work by some of the authors show that this is not a requirement if a modification to MLMC is made \cite{Det18}. The adaptation, called Continuous Level Monte Carlos (CLMC), allows sample-dependent adaptive grids to be built; allowing the multilevel framework to exploit the computational advantages of adaptive finite elements. There is significant future opportunities in stochastic composite analysis to exploit this extension, since often defects arising in manufacturing are localised.

\subsection{Implementation of MLMC}

In this section we discuss how the MLMC algorithm can be implemented in practice, and how the (optimal) values of $L$, $M_\ell$ and $\{N_\ell\}_{\ell=0}^L$ can be computed `on the fly' from the sample averages and the sample variances of $Y_\ell$. For ease of presentation, let us define $Y_0 = Q_{M_0}$. We will also restrict ourselves to the case of uniform mesh refinement where the mesh size is simply halved each time, i.e. $h_\ell = 2^{-\ell} h_0$, but this is not necessarily required \cite{Det18}.

\smallskip
We wish to estimate $\mathbb E[Q]$ within a prescribed RMSE $e$, which is made up of two parts, the bias error and the sampling error \eqref{eq:MC_MSE_2}. Firstly to estimate the bias error, let us assume that $M_\ell$ is sufficiently large, so that we are in the asymptotic regime so that
\begin{equation}
\Bigl|\mathbb{E}[Q_{M_\ell}-Q] \Bigr| \sim M_\ell^{-\alpha}.
\end{equation} 
The number of degrees of freedom on level $\ell$ is given by $M_\ell \approx m^\ell M_0$. For the two-dimensional numerical examples which follow below, domains are discretized by quadrilateral elements so we take $m=4$.  It follows by the reverse triangle inequality that
\smallskip
\begin{equation}
\Bigl|\mathbb E[Y_\ell]\Bigr| = \Bigl|\mathbb E[Q_\ell - Q_{\ell-1}]\Bigr| = \Bigl|\mathbb E[Q_{\ell-1} - Q] - \mathbb E[Q_{\ell} - Q]\Bigr| \geq \Bigl| |\mathbb E[Q_{\ell-1} - Q]\Bigr| - \Bigl|\mathbb E[Q_{\ell} - Q]| \Bigr|.
\end{equation}
By noting that $\mathbb E[Q_{l-1} - Q] \geq cm^\alpha \;\mathbb E[Q_l - Q]$, for some constant $c \approx 1$, we get
\begin{align}
|\mathbb E[Y_\ell]| \geq  (cm^\alpha - 1)\Bigl|\mathbb E[Q_\ell - Q]\Bigr|
\end{align}
Rearranging this expression for the bias error $|\mathbb{E}[Q_{M_\ell}-Q]|$, setting $c=1$ and approximating $|\mathbb{E}[Y_\ell]|$ by the Monte Carlo estimate $\hat Y^{MC}_{\ell,N_\ell}$, the bias error on level $\ell$ can be over-estimated by
\begin{equation}\label{eqn:bias_estimate}
\Bigl|\mathbb{E}[Q_{M_\ell}-Q] \Bigr|  \le \frac{1}{cm^\alpha - 1}\Bigl|\mathbb E[Y_\ell]\Bigr| \le \frac{1}{m^{\alpha}-1} \widehat{Y}^{\text{MC}}_{\ell,N_\ell}\,.
\end{equation}
The sample variance is estimated in the standard way
\begin{equation}\label{eqn:varianceestim}
s^2_\ell = \left(\frac{1}{N_\ell}\sum_{j=1}^{N_\ell} (Y_\ell^{(j)})^2\right) - \left(\widehat{Y}^{\text{MC}}_{\ell,N_\ell}\right)^2 \approx V_\ell\,.
\end{equation}
We summarise the adaptive method in Algorithm 1 and note note that since each sample is independent, Algorithm 1 can be readily parallelized by distributing samples across all processors.

\begin{algorithm}
  \caption{Multilevel Monte Carlo Algorithm}
\label{alg:mlmc}
  \begin{algorithmic}[1]
    \STATE Set $e,\theta,N^\star$
      \STATE Set $L = -1$ \& \texttt{converged = false}
      \WHILE{\texttt{converged == false}}
        	\STATE Compute $N_L = N^\star$ samples on level $L$.
			\STATE  Estimate $V_\ell$ from samples on levels $\ell$, using \eqref{eqn:varianceestim}.
            \FOR{$\ell = 0$ to $L$}
        \STATE {\bf Estimate} optimal samples $\hat{N}_{\ell}$ on level $\ell$ using \eqref{eqn:optimal_Nl}.
        \IF{$N_\ell < \hat{N}_\ell$}
         \STATE Compute $\hat{N}_\ell - N_\ell$ more samples on level $\ell$.
         \STATE Set $N_\ell = \hat{N}_\ell$.
         \ENDIF
      \ENDFOR
      	 \STATE Estimate bias $\hat{e}_b$ on level $L$ using \eqref{eqn:bias_estimate}.
	 	 \IF{$\hat{e}_b < e_b$}
		 	\STATE Set \texttt{converged = true}
		\ENDIF
      \ENDWHILE
  \end{algorithmic}
\end{algorithm}

\subsection{Multilevel Monte Carlo Simulation with Selective Refinement for the Computation of Failure Probabilities}\label{sec:MLMC_SR}

\smallskip\noindent
For many engineering applications, estimating the expected value of a specific quantity is of limited interest, instead often we wish to compute the probability that the failure load $\lambda$ is less than a `safe' load $\lambda^\star$. Within the MLMC framework, the quantity of interest is then the binomially distributed random variable $Q = \boldsymbol{1}(\lambda < \lambda^\star)$, which takes value $1$ if $\lambda < \lambda^\star$ and $0$ otherwise. The failure probability can then be approximated by evaluating $\mathbb E[Q] = \mathbb P(\lambda < \lambda^\star)$. 

\smallskip\noindent
For aerospace applications these probabilities are necessarily small, and obtaining good estimates for these rare events is difficult since, by definition, a large number of samples are required to observe even a single case. One of the main issues is that a simple binomial distribution ($Q = 0$ or $1$) loses important information regarding how close a given sample is to failing; in particular, $Q(\lambda)$ is a step function at $\lambda^\star$. One proposed method for improving convergence is to use a smooth quantity of interest, which takes intermediate values between $0$ and $1$ if $\lambda$ is close to the critical value \cite{Gil08}. Here, however, we propose a different approach which combines the error estimator in Eqn.\,\eqref{eqn:bias_estimate} and the multilevel framework presented in Section \ref{subsec:MLMC}, motivated from an approach proposed by Elferson et al. \cite{Elf16}.

\begin{figure}\label{fig:failureProb}
\centering
\includegraphics[width=0.7\linewidth]{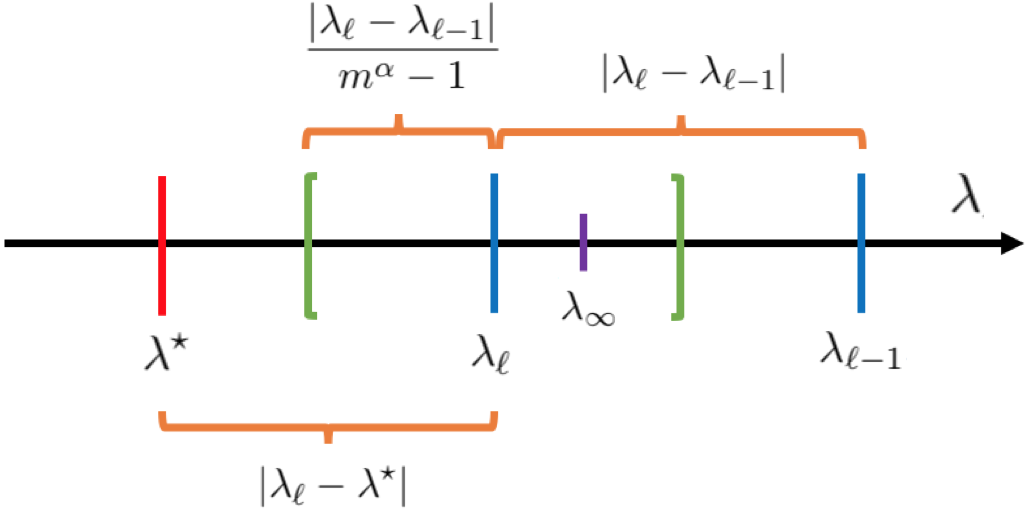}
\caption{Graphical representation of MLMC with selective refinement}
\end{figure}

\smallskip\noindent
Following a similar calculation to \eqref{eqn:bias_estimate}, it follows that the bias error for a given sample can be estimated by
\begin{equation}\label{eqn:sampleBias}
|\lambda_{\ell} - \lambda| \approx  \frac{|\lambda_{\ell} - \lambda_{\ell - 1}|}{m^\alpha - 1}\; .
\end{equation} 
Therefore, if we wish to approximate $Q_\ell$ but observe that for some level $i < \ell$  
\begin{equation}
|\lambda_i - \lambda^\star| \geq \frac{|\lambda_{i} - \lambda_{i-1}|}{m^\alpha - 1}\; ,
\label{eqn:testFailure}
\end{equation}
then Eqn.\,\eqref{eqn:sampleBias} ensures that $Q_\ell = Q_i$ for all $\ell \geq i$. In many cases,  it is then unnecessary to calculate $\lambda_\ell$ on high levels in order to obtain the `fine' level approximation of $Q_\ell$, since the coarse approximation is sufficiently far from  $\lambda^\star$ (as illustrated in Fig.\,\ref{fig:failureProb}). This selective refinement technique is summarised in Algorithm \ref{alg:mlmc_sr}. The key point is that this modification simply reduces the average cost per sample on refinement levels $\ell > 2$, whilst the original multilevel algorithm (as described by Algorithm 1) remains unchanged.

\smallskip
Elferson et al. \cite{Elf16} first formalised the gains of MLMC-SR over MLMC and std. MC for computing failure probabilities, showing that the expected cost to compute one level $\ell$ sample of the failure probability functional
$Q_\ell$ using the selective refinement method (as presented in Algorithm 2) can be bounded by 
\begin{equation}\label{eqn:costpersample}
\text{Cost}(Q_\ell) \leq C_5 M_\ell^{\gamma-\alpha}
\end{equation}
The order of growth, with respect to degrees of freedom $M_\ell$, is shown to be significantly reduced from the corresponding cost per sample for the MLMC simulation. This is due to the fact that only a fraction of level $\ell > 0$ samples are solved on their highest refinement levels, with work instead concentrated on the lower (computationally cheaper) levels.

\begin{algorithm}
\caption{Selective refinement procedure for one sample of a failure probability calculation.}
\label{alg:mlmc_sr}
  \begin{algorithmic}[1]
\STATE For given $\ell$, $i$ and $\lambda^\star$
\FOR {levels $j = 0,\ldots,\ell$}
 \STATE Compute $\lambda^{(i)}$ on level $j$.
 \IF {$ j > 1$} 
 	\IF{$|\lambda_j^{(i)} - \lambda^\star| \geq |\lambda^{(i)}_j - \lambda^{(i)}_{j-1}| / (m^\alpha - 1)$}
 		\STATE Set $\lambda^{(i)}_\ell = \lambda_j^{(i)}$.
		\STATE Exit {\bf for loop}
 	\ENDIF
 \ENDIF
\ENDFOR
\STATE Evaluate failure probability functional $Q^{(i)}_\ell = \boldsymbol{1}(\lambda^{(i)}_\ell - \lambda^\star)$
\end{algorithmic}
\end{algorithm}

\smallskip\noindent
We note that the standard estimates for the mean and variance of $Y_\ell$ deteriorate as $\ell$ increases.  This is a significant practical challenge that arises when computing failure probabilities for both MLMC and MLMC-SR algorithms, since as $\ell$ increases the probability of $Y_\ell \neq 0$ approaches zero. We are particular interested in finding stable estimates for (very) small probabilities. It is important that they are not underestimated, since they are used to bound the numerical bias and sampling error which control the stopping criterion for the MLMC algorithm. To address this challenge we directly recap the ideas introduced by Elferson {\em et. al} \cite{Elf16}.

\smallskip\noindent
In general, $Y_\ell$ is a trinomial random variable taking values of either $-1$, $0$ and $1$. If $p_{+1}$ is the probability observing a failure of level $\ell$ and not on $\ell-1$ (i.e. $Y_\ell=1$) and $p_{-1}$ probability of failure on $\ell-1$ but not on $\ell$ (i.e. $Y_\ell=-1$); then $p_{+1}$ and $p_{-1} \rightarrow 0$ as $\ell \rightarrow \infty$. Therefore the accuracy of sample estimates for the mean and variance of $Y_\ell$ deteriorates. The true values are
\begin{equation}
\mathbb E[Y_\ell] = p_{+1} + p_{-1}  \quad \text{and} \quad \mathbb V[Y_\ell] = p_{+1} + p_{-1} + (p_{+1}-p_{-1})^2,
\end{equation}
In \cite{Elf16} the following biased estimators $\tilde{p}_{+1}$ for the parameter $p_{+1}$ is introduced to overcome this issue
\begin{equation}\label{eqn:biased_estimator}
\tilde p_{+1} = \frac{x_{+1} + k}{N_\ell+k}
\end{equation}
where $x_{+1}$ denotes the number of samples for which $Y_\ell = +1$ within $N_\ell$ samples and $k \in \mathbb N$. An identical expression is used for $\tilde p_{-1}$. To quantify the accuracy of these estimators we calculate the relative variance  $\mathbb V[\tilde{p}]/\mathbb{E}[\tilde p]^{2}$, for which a value greater or equal to one indicates a significant departure from the trinomial distribution. We see that for this choice of biased estimator \eqref{eqn:biased_estimator}, the value is less than one:
\begin{equation}
\frac{\mathbb V[\tilde p_{+1}]}{\mathbb E[\tilde p_{+1}]^2} = \frac{N_\ell p_{+1}(1-p_{+1})}{(N_\ell p_{+1}+k)^2}\leq \frac{N_\ell p_{+1}}{(N_\ell p_{+1}+k)^2}<1
\end{equation}
Choosing a large value of $k$ gives a large bias in the estimator, but a smaller relative variance. The bias of the estimator is significant if  $N_\ell p_{+1} \ll k$ and there are too few samples to estimate $p_{+1}$ accurately. However, $\tilde{p}_{+1}$ acts as a bound in that case.

\smallskip
\noindent{\em Remark:} For the special case of a stochastic eigenvalue problem with a nested hierarchy of grids (considered in Sec. \ref{sec:results}) the Min-max Principle \cite{Parlett} ensures that $\lambda_{\ell} \leq \lambda_{\ell-1}$. It naturally follows that $p_{-1} = 0$ and $Y_\ell$ is a binomial random variable.

\section{Test Problem I : Compressive Strength of Fibre Composites with Random Fibre Misalignment}

\begin{figure}[!h]
\centering
\includegraphics[width = 0.45\linewidth]{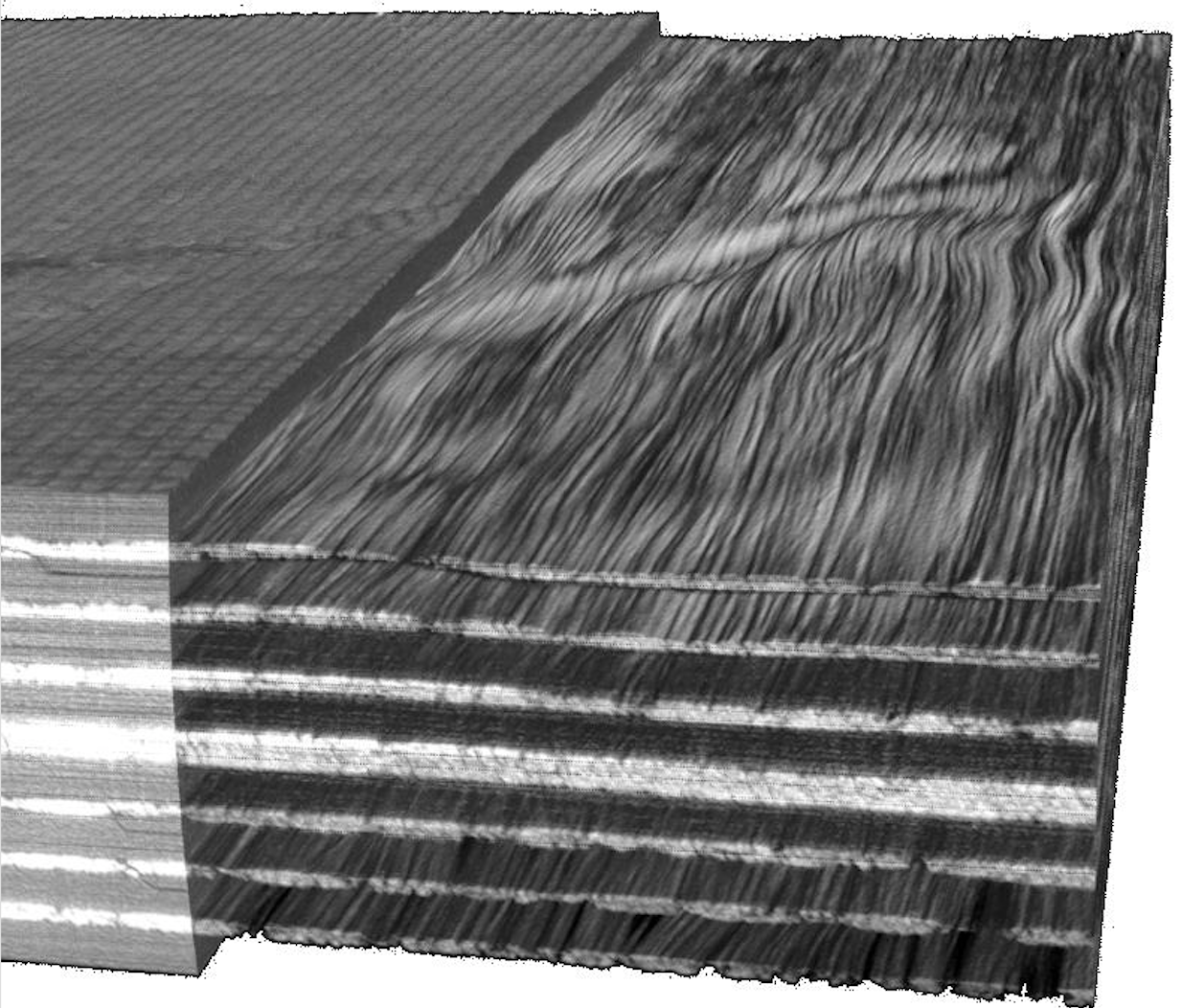}\qquad \includegraphics[width = 0.45\linewidth]{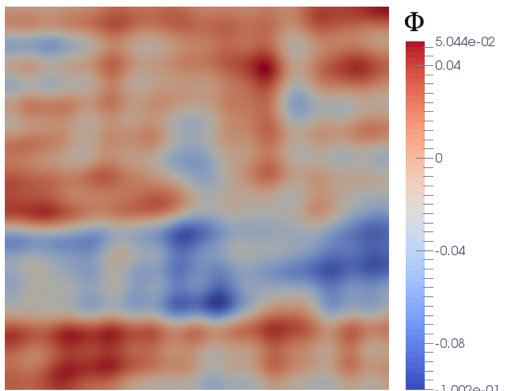}
\label{fig:waviness}
\caption{(Left) CT image showing random fibre waviness within a composite laminate. (Right) Sample of the random waviness field $\Phi$ with $N_{KL} = 400$ and covariance parameters as taken in the results section \eqref{eqn:covarianceParameters}.}
\end{figure}

It is well established that the compressive failure of undamaged composites is primarily governed by plastic micro-buckling  (or kinking) of the fibres \cite{Bud83,Liu04}, and this failure is initiated in regions of local fibre misalignment or waviness. The classical micro-mechanical model for the compressive strength of a composite $\sigma$ given by Budansky \cite{Bud83} is \begin{equation}
\sigma = \frac{G}{1 + |\Phi|/\gamma_y}
\label{eqn:BudModel}
\end{equation}
where $G$ , $\Phi$ and $\gamma_y$ are the shear modulus, fibre misalignment angle and shear strain at failure, respectively. This idealised model (often referred to as \textit{kinking theory}) assumes the misalignment or kink of known angle $\Phi$. Observations of real fibre waviness show that the misalignment $\Phi$ is not a single value, but a complex random field, as seen in Fig. \ref{fig:waviness} (left). In practical applications it is then unclear what value of misalignment $\Phi$ should be used in \eqref{eqn:BudModel}; possible options include the root mean square or the maximum misalignment. In fact, the compressive strength is also a random variable, with a distribution intricately coupled with the statistical distribution of $\Phi$. We model the uncertainty in the angle with a spatial random field, as shown for example in Fig. \ref{fig:waviness} (right).

\smallskip
\noindent
In this section we model the uncertainty in the angle via a spatial random field parameterised by observed statistics \cite{Sut12}, and demonstrate the computational savings of the MLMC by estimating the expected compressive strength of a composite $\mathbb E[\sigma]$, with random fibre waviness .


\subsection{A Two-Dimensional Cosserat Continuum Model for a Composite with Random Misalignment}\label{subsec:cosserat}
\noindent
In this test problem we consider a square domain $\Omega$ in the $(x,y)$ plane, made up of uni-directional composite pre-preg. Individual fibres are misaligned by an angle $\Phi(\textbf x)$ to the $x$-axis. This misalignment is modelled as a random field on $\Omega$. The mean and the covariance structure of $\Phi$ will be inferred from measurements of fibre misalignment of carbon fibre pre-pregs available in the literature \cite{Sut12}. The random field $\Phi$ is characterised by a two-point exponential covariance function 
\begin{equation}
k(\textbf x, \textbf y) = s_{\Phi}^2 \exp \left( -\frac{|x_1 - y_1|}{\omega_1} -\frac{|x_2 - y_2|}{\omega_2}\right).
\label{eqn:covarianceFunction}
\end{equation}
The parameters $s_\Phi^2$ and $\omega_i$ denote the variance and correlation length (in each direction) of the misalignment field. To generate a single realisation of this random field, we represent the random misalignment fields as a set of random variables using a Karhunen-Loeve (KL) expansion, an expansion in terms of a countable set of basis functions $\phi_n(x)$ parameterised by independent standard Gaussian random variables $\{\xi_n\}_{n \in \mathbb N}$, given by
\begin{equation}
\Phi(\textbf x) = \sum_{n=1}^{\infty}\sqrt{\mu_n}\phi_n(\textbf x)\xi_n.
\label{eqn:inftyKL}
\end{equation}
Here $\{\mu_n\}_{n \in \mathbb N}$ and $\{\phi_n\}_{n \in \mathbb N}$ are the eigenvalues and associated (normalised) eigenfunctions of the covariance function \eqref{eqn:covarianceFunction}. We note that the eigenvalues $\{\mu_n\}_{n \in \mathbb N}$ are positive and strictly decreasing which provides a natural ordering of the importance of the contribution of each term to $\Phi(\textbf x)$. In a computational setting it is therefore natural to truncate the KL-expansion \label{eqn:inftyKL} after $N_{KL}$ terms, giving a parameterisation of the random field by the set of variables $\boldsymbol \xi = [\xi_1,\xi_2,\ldots,\xi_{N_{KL}}]$. Figure \ref{fig:waviness} (right) shows a realisation of the random field generate using the approach described. For further details of random fields and their implementation within structural applications we refer the reader to the classical text by Spanos and Ghanen (2003) \cite{Spa03}. We note that it is possible to implement more complex covariance functions, and implement them on more complex geometries, see for example \cite{Sca18}.
%

\smallskip
A key consideration when modelling the mechanics of such a composite, is that the shear stiffness parallel to the fibres is an order of magnitude less than the shear stiffness orthogonal to them; and hence, in general, the stress state is non-symmetric, i.e. $\sigma_{12} \neq \sigma_{21}$. As a result, a finite size element of composite carries a coupled stress (a moment per unit area), and the fibres bend to achieve moment equilibrium. A classical approach to capturing these internal bending effects is to model the composite as a Cosserat Continuum \cite{Liu04,Dod15b}. Here, under plane-strain assumptions, each material point has the conventional displacement degrees of freedom $u_1$ and  $u_2$ ($u$ and $v$ in global coordinates), as well as an independent (Cosserat) rotational degree of freedom $\theta_3$. Under the assumption of small deformations and rotations, this gives the small Cosserat strain and curvature measures
\begin{equation}
\varepsilon_{ij} = \frac{du_i}{dx_j} + e_{ijk}\theta_k \quad \mbox{and} \quad \kappa_{ij} = \frac{d\theta_3}{dx_j},
\end{equation}
where $e_{ijk}$ denotes the permutation tensor. The permutation tensor is defined as $e_{123} = e_{312} = e_{231} = 1$, $e_{213} = e_{132} = e_{321} = -1$ and $e_{ijk} = 0$ if any indices are repeated, e.g. $e_{112} = 0$. These strain and curvature measures are work conjugates to the Cosserat stresses $\sigma_{ij}$ and coupled-stress $m_{ij}$, respectively. We introduce  the linear Cosserat constitutive relationships (derived in \cite[Sec. 2.3]{Dod15}), which are expressed in matrix form as $$\boldsymbol \sigma = \textbf C \boldsymbol \varepsilon \quad \mbox{and} \quad \textbf m = \textbf D \boldsymbol \kappa.$$ The matrices $\textbf C$ and $\textbf D$ can be rotated by the misalignment angle $\Phi(\textbf x)$ to the global $x$-axis via the transformation matrices
$ \textbf T^{\varepsilon}_{\phi(\textbf x)}$ and 
$\textbf T^{\kappa}_{\phi(\textbf x)}$
so that the global matrices become
\begin{equation}
\textbf C^*_{\phi(\textbf x)} = (\textbf T^{\varepsilon}_{\phi(\textbf x)})^{-1}\; \textbf C \;\textbf T^{\varepsilon}_{\phi(\textbf x)} \quad \mbox{and} \quad \textbf D^*_{\phi(\textbf x)} = (\textbf T^{\kappa}_{\phi(\textbf x)})^{-1}\; \textbf D \;\textbf T^{\kappa}_{\phi(\textbf x)}.
\end{equation}
The force and moment equilibrium equations for a small element of composite, in the absence of body forces and coupling are given by
\begin{equation}
\frac{d\sigma_{ij}}{dx_j} = 0 \quad \mbox{and} \quad \frac{dm_{ij}}{dx_j} + e_{ijk}\sigma_{jk} = 0.
\label{eqn:equilbrium}
\end{equation}
In our model, these equilibrium equations are subject to the Dirichlet boundary conditions
\begin{equation}
u(\textbf x) = 0 \; \mbox{on} \; x = 0  \quad \mbox{and} \quad  u(\textbf x) = \Delta < 0 \; \mbox{on} \; x = L,
\end{equation}
and
\begin{equation}
v(\textbf x) = 0 \: \text{ on $y = 0$ and $y = L$.}
\end{equation}
To solve \eqref{eqn:equilbrium} using the finite element method, the differential equations are recast as a variational problem. We seek a solution $(\textbf u, \theta_3) \in V^2 \times W$, such that for all test functions $({\textbf v},\hat \vartheta_3) \in V^2 \times W$ the equality
\begin{align}\label{eqn:VarCosserat}
\int_\Omega \textbf C^*_{\phi(\textbf{x})}\varepsilon(\textbf u,\theta_3):\varepsilon(\textbf v,\vartheta_3) + \textbf D^*_{\phi(\textbf{x})}\;\kappa(\theta_3):\kappa(\vartheta_3)\;d\textbf x
= \int_\Gamma \textbf t \cdot \textbf v + \mu\vartheta_3\;d\textbf x,
\end{align}
holds. Here $\textbf t$ denote the stress traction, and $\mu$ the coupled stress traction on the boundary of the domain $\Gamma$. The spaces $V$ and $W$ are  appropriate function spaces on which the components of $\textbf u$ and the Cosserat rotation $\theta_3$ are defined. Here, an appropriate choice is the Sobolev Space $H^1$; that is, the space of all square integrable functions with square integrable first derivatives satisfying the boundary conditions.

\smallskip\noindent
To approximate \eqref{eqn:VarCosserat}, the domain $\Omega$ is uniformly discretized into a set of 4-node quadrilateral elements $$\mathcal Q_h = \{\Omega_e^{(i)}\}_{i=1}^{\rm{nel}},$$ where $\rm{nel}$ denotes the number of elements and $h$ is the side-length of the elements. The solution is approximated by restricting \eqref{eqn:VarCosserat} to the finite dimensional subspace $V_h^2 \times W_h \subset V^2 \times W$. In these examples $V_h$ and $W_h$ are chosen to be the set of piecewise bi-linear functions on $\mathcal Q_h$, and we denote the corresponding approximate solution by $\textbf u_h$ and $\theta_{h}$. As for any standard finite element analysis, substitution of the approximations $\textbf u_h$ and $\theta_h$ allows \eqref{eqn:VarCosserat} to be rewritten as a linear system of the form 
\begin{equation}
\textbf K \textbf d = \textbf f
\label{eqn:matrixinplane}
\end{equation}
where $\textbf K \in \mathbb{R}^{M\times M}$ is the global stiffness matrix and $\mathbf f \in \mathbb{R}^M$ is the load vector due to the prescribed boundary conditions. The vector $\textbf d \in \mathbb{R}^M$ contains the coefficients of all degrees of freedom in the expansions of $\textbf u_h$ and $\theta_h$ above. If $\texttt{nnode}$ is the total number of nodes in the grid, then $M = 3 \texttt{nnode}$.

\smallskip
From the solution $\textbf d$, we wish to calculate the compressive strength, $\sigma$, of the composite. We consider the quadratic failure criterion \cite{Bud83}, and introduce the effective stress $\tau_e$, which is defined in terms of the transverse stress $\sigma_{22}$ and the shear-stress parallel to the fibres $\sigma_{12}$. In particular,
\begin{equation}
\tau_e = \sqrt{\sigma_{12}^2 + \left(\frac{\sigma_{22}}{R}\right)^2}.
\end{equation}
The material parameter $R$ is the ratio of the transverse and shear yield strength of the material. We say that failure occurs when the effective stress is equal to the shear strength of the material $\tau_y$; i.e. $\tau_e = \tau_y$.

\smallskip
In the results which follow, we estimate the compressive strength by first computing $\mathbf{d}$ for a prescribed compressive end-shortening $\Delta$. In order to remove the influence of boundary conditions, we then find the maximum value $f^*$ of $f = \tau_e/\tau_y$ over all integration points within elements contained in a central square subregion $\Omega'$ of $\Omega$, which has area $|\Omega'|$. As the problem under consideration is linear, the compressive strength $\sigma$ is then given by
\begin{equation}
\sigma = \frac{f^*}{|\Omega'|}\int_{\Omega'}\sigma_{x}\;d\textbf x.
\label{eqn:failureCriterion}
\end{equation}

\subsection{Results}

For the experiments that follow, we consider material parameters for unidirectional pre-preg AS4/8552, with material constants taken from the Hexcel Data Sheet \cite{Hex08}. Specifically, we take 
$$v_f = 0.59, \quad E_f = 230\mbox{GPa}, \quad E_m=9.25\mbox{GPa}, \quad G_f = 95.83\mbox{GPa},$$
$$G_m = 5.13\mbox{GPa}, \quad d = 7\mu\mbox{m} \quad \mbox{and} \quad \tau_c = 114\mbox{MPa}.$$
The stochastic model for random misalignment is parameterised based on data in the literature; in particular, the measurements of in-plane waviness in pre-preg given by Sutcliffe et al. \cite{Sut12} which agree well with other values given by \cite{Liu04,Jel92,Wis94}. In this paper the correlation lengths, $\omega_1$ and $\omega_2$, are defined differently to those given in Sutcliffe et al. \cite{Sut12}, which we will denote by $\omega^\star_1$ and $\omega^\star_1$. They defined the correlations lengths as the lag at which the auto-correlation function is equal to $0.1$, i.e. when $k(\textbf x,\textbf y) /\sigma_\phi^2 = 0.1$, and therefore $\omega_1= -\omega_1^\star / \log(0.1)$ and similarly for $\omega_2$. Therefore the covariance function \eqref{eqn:covarianceFunction} is parameterised with the following values
\begin{equation}
\omega_1 =  229d ,\quad \omega_2 = 61d \quad \mbox{and} , \quad \mbox{and} \quad s_\Phi = 0.035\rm{rad} 
\label{eqn:covarianceParameters}
\end{equation}
Figure \ref{fig:waviness} shows a random field generated with the above parameters. Having fixed the correlation lengths of the wrinkles, the domain size is chosen to be $L = 2.5\omega_1$. Furthermore, $|\Omega'|$ (as introduced in \eqref{eqn:failureCriterion}) is chosen to be the square subdomain centred in $\Omega$ with sides of length $1.25\omega_1$. The coarsest finite element grid (level $\ell = 0$) has a mesh size of $h_0 = L/8$ (i.e. with 64 elements and $M_0 = 243$ degrees of freedom), and subsequent grids are generated by uniform refinement as shown in Fig.\,\ref{fig:meshes}. The number of KL modes is also increased with the levels $\ell$, $N_{KL}^\ell = 50 + 50 \ell$.

\smallskip
Before comparing the MLMC algorithm with standard MC, we first estimate how the computational cost scales will $M_\ell$ to estimate the parameter $\gamma$. By recording times to compute 100 samples from level $\ell = 0$ to $\ell = 5$; from this we estimate 
\begin{equation}
\mathcal C_\ell \leq CM_{\ell}^{1.3},
\end{equation}
where $C$ is some constant independent of $M_{\ell}$, and hence we take $\gamma= 1.3$ in our calculations to follow.

\begin{figure}
\centering
\includegraphics[width=\linewidth]{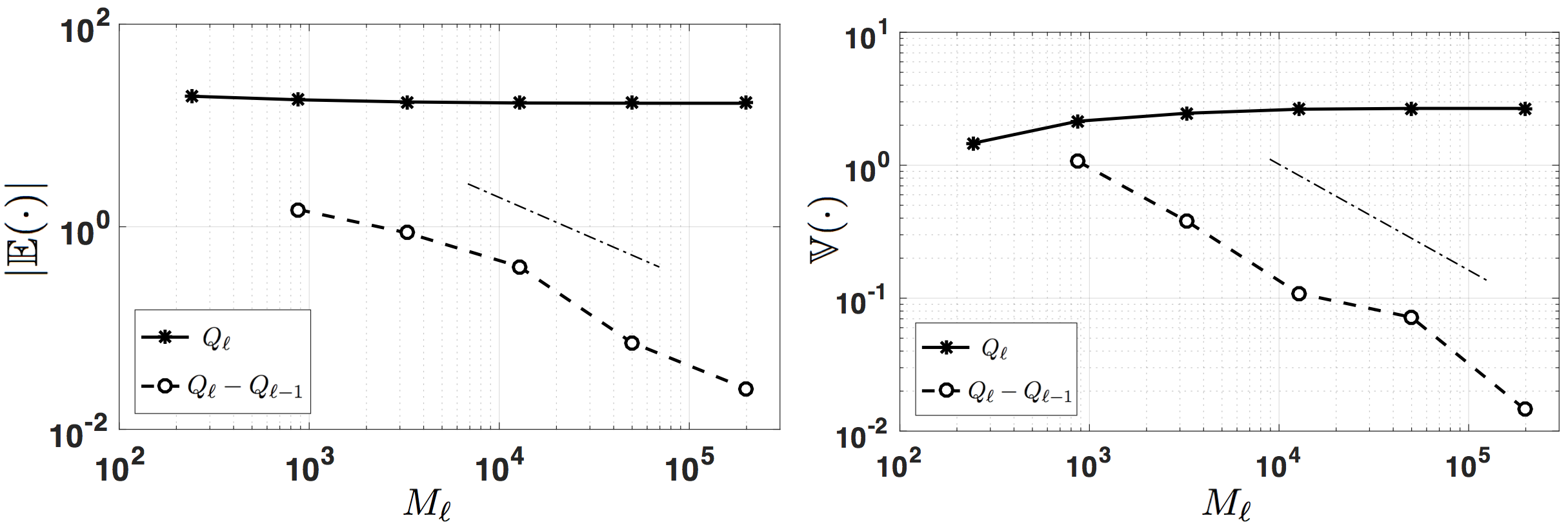}
\caption{(Left) Expected value of $Q_\ell$ and $Y_\ell = Q_\ell - Q_{\ell -1}$ against degrees of freedom $M_\ell$, dashed line shows $\alpha \approx 0.786$. (Right) Variance of $Q_\ell$ and $Y_\ell = Q_\ell - Q_{\ell -1}$ against degrees of freedom $M$, dashed line shows $\beta \approx 0.740$.}
\label{fig:test1}
\end{figure}

\smallskip\noindent
We carry out a series of MLMC simulation, where we take our Quantity of Interest as
$Q = \sigma$ as defined by \eqref{eqn:failureCriterion} over a range of tolerance values from $e=42.13$MPa down to $3.06$MPa (relative error of $3\%$ down to $0.2\%$). The values of the parameters $\alpha$ and $\beta$, as defined in Section \ref{subsec:MLMC}, can be determined from Fig.\,\ref{fig:test1}, which shows the log-log plots of the mean and variance of $Q_\ell$ and $Y_\ell = Q_\ell - Q_{\ell-1}$, with respect to the total number of degrees of freedom, $M_\ell$. Looking first at the behaviour of the expectation of $Q_\ell$ and $Y_\ell$ (left), we see that
\begin{equation*}
\mathbb{E}[Y_\ell] \leq CM_\ell^{-0.786}
\end{equation*}
approximately, and hence $\alpha \approx 0.786$. Next, considering the variance plot (centre), we see that
\begin{equation*}
\mathbb{V}[Y_\ell] \leq CM_\ell^{-0.740}
\end{equation*}
approximately, and hence $\beta \approx 0.740$.

\begin{figure}
\centering
\includegraphics[width = 0.55\linewidth]{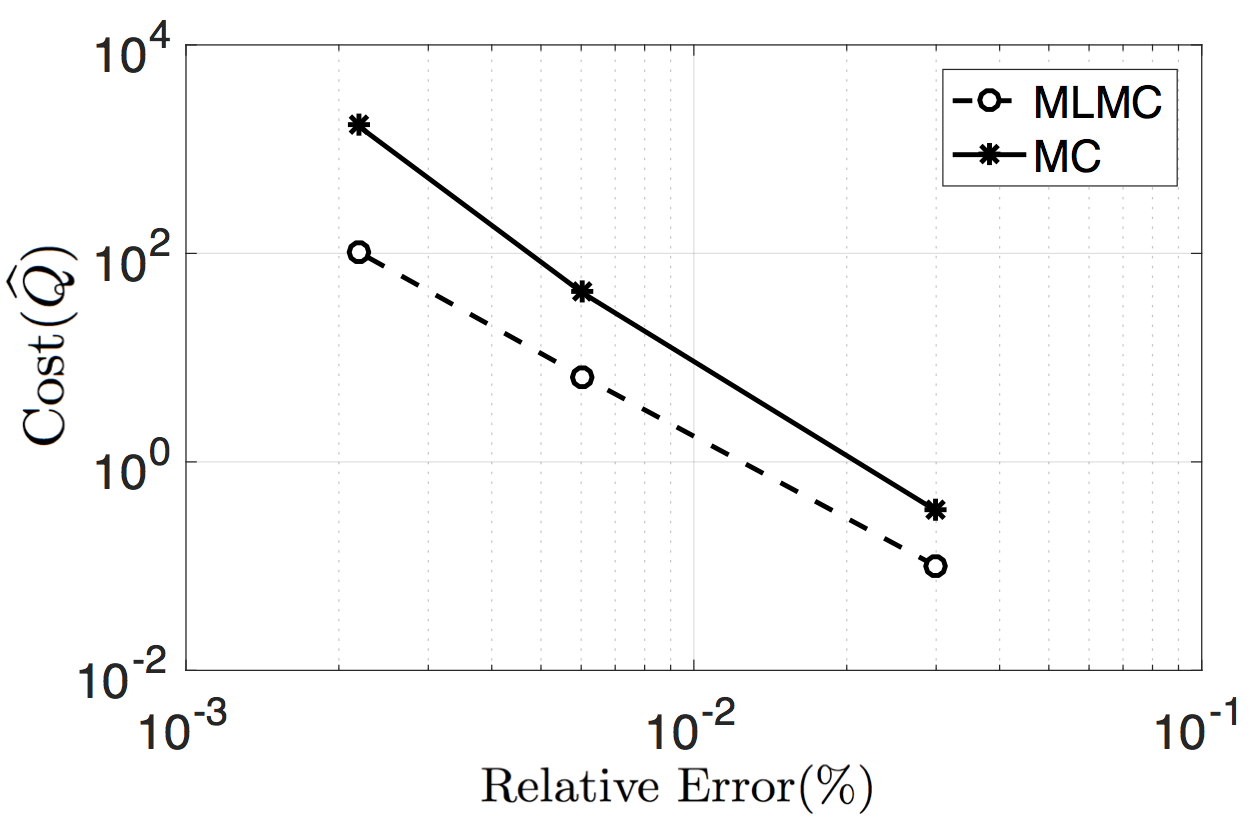}\includegraphics[width = 0.4\linewidth]{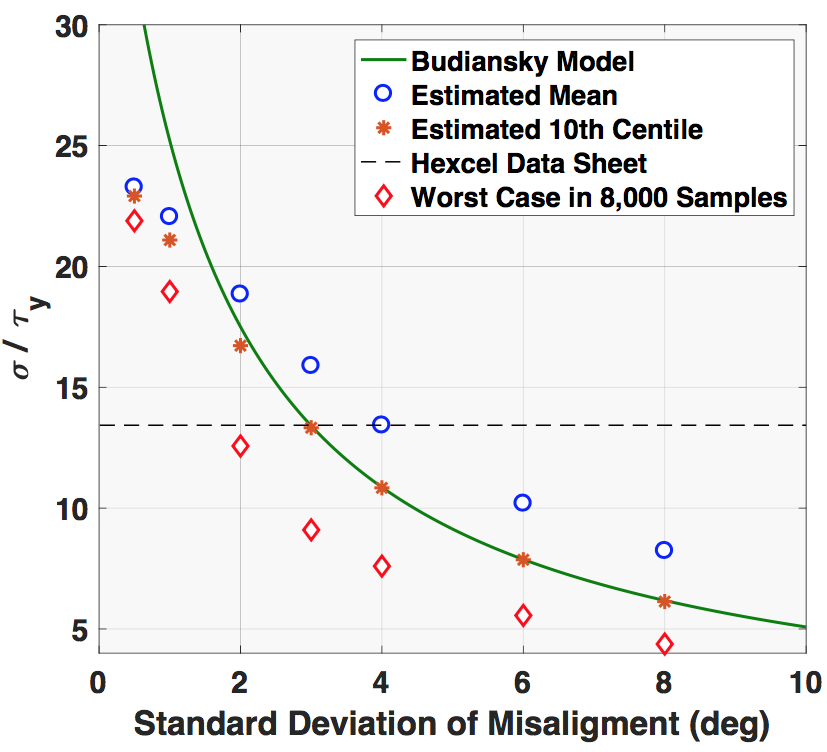}
\caption{(Left) Relative error ($\%$) against computational cost for standard MC (Cost $\sim e^{-3.22}$) and MLMC simulations (Cost $\sim e^{-2.64}$) (Right) Normalised Compressive Strength against standard deviation of the misalignment field. 
}
\label{fig:main_result}
\end{figure}

\smallskip
Figure \,\ref{fig:main_result} (Left) compares the computational cost of the MLMC simulation versus standard MC, with respect to error tolerance. For the current parameter values, Eqn.\,\eqref{eqn:mlmcComplexity} predicts the cost of the MLMC simulation to grow proportionally to $e^{-2.68}$, whilst the cost of the standard MC simulation grows like $e^{-3.64}$. The numerical experiment verifies these predictions; considering the gradients of the plots, we see that the cost of the MLMC simulation is approximately proportional to $e^{-2.64}$, whilst that of the MC simulation is proportional to $e^{-3.22}$. These gains are explicitly quantified in Table \ref{tab:test1costs}, which lists the optimal numbers $N_{\ell}$ of samples on each refinement level for three absolute error tolerances, as given by \eqref{eqn:optimal_Nl}, along with the total computational costs of the MLMC simulation. Also given are the required numbers of samples for the corresponding standard Monte Carlo simulations, from which a computational speed-up factor may be calculated. In particular, we see that for an absolute error of $e = 3.06$MPa the MLMC algorithm reduces the computational cost by a factor of $16$ over standard MC; in absolute terms this reduces computation times from $28$ hours to under $2$ hours.

\begin{table*}
\centering
\begin{tabular}{|c|c|c|c|c|c|c|c|c|c|}
\hline
\multirow{2}{*}{\begin{tabular}[c]{@{}c@{}}$e$ \end{tabular}} & \multicolumn{6}{c|}{$N_\ell$}        & \multirow{2}{*}{\begin{tabular}[c]{@{}c@{}}MLMC \\ Cost \end{tabular}} & \multirow{2}{*}{\begin{tabular}[c]{@{}c@{}}MC \\ Cost \end{tabular}} &  \multirow{2}{*}{\begin{tabular}[c]{@{}c@{}}Saving \\ Factor \end{tabular}} \\ \cline{2-7} & 0 & 1 & 2 & 3 & 4 & 5 &  &    &  \\ \hline
3.01\% & 513 & 237 & 34 & 8 & - & - & 0.10   & 0.34  & 3.32 \\
0.63\% & 22,014 & 6,191 & 1,449  & 337  & 123  & - & 6.65    & 42.26 & 6.36\\
0.22\% & 240,427 & 67,611 & 15,822 & 3,684 & 1347 & 283 & 103.84  & 1685.50 & 16.23 \\ \hline
\end{tabular}
\caption{Cost comparison between std. MC and MLMC for test I.}
\label{tab:test1costs}
\end{table*}

\smallskip
Whilst the principle aim of this paper has been to demonstrate the computational savings of the MLMC methodology, we now also compare the results to theoretical and experimental work in the literature. Firstly we consider the influence of the size of the standard derivation $s_\Phi$ of the misalignment field on the compressive strength of  AS4 /8552. Using the new multilevel methodology (with $L=4$), Fig. \ref{fig:main_result} shows the estimated mean $\mathbb E(\sigma)$ (blue markers), estimated 10th percentile (red markers) strength values and the worst case in 8,000 samples on level $4$ with an increasing standard deviations of the misalignment field $s_\Phi$. The results are compared to the classical Budiansky `kinking' model \eqref{eqn:BudModel} and also the Hexcel data sheet value for AS4/8552 ($\sigma/\tau_y =13.43$). The estimates for the 10th per centile for $s_\Phi \geq 2$ agree very well the Budiansky model,  where both predict a significant decrease in compressive strength with increasing fibre misalignment. Discrepancies between strength values at lower misalignment angles, suggest a small misalignment angles is not dominated by shear, but by failure in the $\sigma_{22}$ direction, which is not accounted for in the Budiansky model~\cite{Bud83}.  

\section{Test Problem II - Buckling performance of a wing skin panel
  with uncertain ply orientations}
\label{sec:model}

In this section we describe a model problem to test the
multilevel Monte Carlo method with selective refinement, as described in Sec. 2 \ref{sec:MLMC_SR}. Here, as an illustrative example for our new methodology,
the structural performance of a wing skin panel subject to a typical
in-service load is considered. Failure of the panel occurs when the panel buckles. Additional different engineering scenerios are given in \cite{But15}.

\subsection{Model Setup and Mathematical Description}

Consider a rectangular composite plate of thickness $t$, length $L_x = 636$mm
and width $L_y = 212$mm, with the un-deformed mid-plane of the plate occupying
the domain $\Omega = [0,L_x] \times [0,L_y]$ with boundary
$\Gamma$. The laminate is made up of $8$ identical, orthotropic, composite
plies characterised by the elastic tensor $\bf{Q}$, thickness $0.8$mm  and arranged in a in a fully uncoupled (Winckler)
stacking sequence $$\boldsymbol \psi = [45^\circ , -45^\circ, -45^\circ, 45^\circ, -45^\circ, 45^\circ, 45^\circ, -45^\circ ].$$ 
The elastic ply properties, are taken from the IM7-8552 data sheet, so that $E_{11} = 130.0$GPa, $E_{22} = 9.25$GPa, $G_{12} = 5.13$GPa, $\nu = 0.36$ and $G = 5.13$GPa.

\smallskip
In this model problem we consider that as-manufactured the ply orientation is uncertain due to angle tolerances in the laying machine. Therefore 
we add a small, constant, random perturbation, $\phi_i$ to each
pristine ply angle $\psi_i$, for $i = 1, \ldots, 8$. In this way, a new
"defective" stacking sequence, $\underline{\psi}^d =
[\psi_1^d,...,\psi_8^d]$, is obtained. We assume that these random angle perturbations
are normally distributed such that $\phi_i \sim \mathcal N(0, 3^2)$. This standard deviation of the perturbations has been chosen to conform with the accuracy of automated fibre placement (AFP) machines in the industry. Typically machines have an allowable error tolerance of $5^\circ$. Hence, in order to obtain
sample perturbations $\phi_i$ satisfying this error tolerance with
$95\%$ confidence, the required standard
deviation is $5^\circ/1.65 = 3^\circ$ (where $1.65 = z_{.05}$ is the
critical $z$ value for the one-sided $95\%$ confidence interval of a
normal distribution).

\smallskip
The deformation of the
plate is described by the
vertical displacement $w(x,y)$ and rotations of the mid-plane
$\underline \theta(x,y) = [\theta,\phi]$. The plate is subjected to
uniform, unit, axial compression stress, whilst being simply-supported around all boundaries. The critical buckling load for the plate is calculated using Reissner-Mindlin (RM) plate theory, since the advantages over Kirchhoff Plate theory are well documented \cite{hughes}. The problem is therefore reduced to a 2D problem in $\Omega$, by applying classical laminate theory (CLT)\cite{gurdal}, which gives the laminate stiffness tensors\begin{equation}
A = \sum_{k=1}^K {\bar Q}^{(k)} (z_k - z_{k-1}), \quad B = \frac
12 \sum_{k=1}^K {\bar Q}^{(k)} (z_k^2 - z_{k-1}^2) \quad \mbox{and} \quad  D = \frac
13\sum_{k=1}^K {\bar Q}^{(k)} (z_k^3 - z_{k-1}^3),
\label{eqn:CLTABD}
\end{equation}
where $z_k$ is the distance from the top edge of the $k^{th}$ ply to the neutral axis of the plate and where ${\bar{Q}}^{(k)}$ is the elastic
tensor of the $k^{th}$ ply in global coordinates. These homogenised
tensors connect in-plane strains $ \varepsilon$ and out-of-plane curvatures $ \kappa(\underline \theta) = \frac{1}{2}\left(\nabla \underline
\theta + \nabla \underline \theta^T\right)$, with in-plane stress and plate bending moments. Under
the additional assumption that the in-plane and out-of-plane behaviour is decoupled, it follows that the in-plane stress 
and the moment are then given by
\begin{equation}
 \sigma = t^{-1}  A \; \varepsilon  \quad \mbox{and}
\quad  \mu =   D^* \; \kappa, 
\end{equation}
respectively. Here, $ D^* =  D -  B^T A^{-1} B$,
which conservatively knocks down the bending resistance of the panel to account for  coupling effects.

\smallskip\noindent
In the absence of body forces, a moment equilibrium for the RM plate gives the linear eigenvalue problem
\begin{equation}
\nabla \cdot ( D^* \kappa(\underline \theta))  - kG(\nabla
w - \underline \theta) = \lambda\nabla \cdot( \sigma\nabla w) \
\quad  \mbox{such that} \quad w = 0\;\; \mbox{and} \;\;  \mu \cdot \underline n = 0 \;\;
\mbox{on}\;\;
\Gamma ,
\label{eqn:momentequil}
\end{equation}
where $G$ is the through thickness shear stiffness and $k = 5/6$ is
the shear correction (both constants), whilst $ \sigma$ is the in-plane stress field  given by 
$$
\sigma = 
\begin{bmatrix}
-1 & 0 \\
0 & 0
\end{bmatrix}
$$ 
Again, \eqref{eqn:momentequil} is solved using FEM, and therefore the weak form of the eigenvalue problem is used, such that the problem becomes: Find the smallest (positive real) eigenvalue $\lambda $ and associated (buckling) eigenmode $0 \neq (\underline \theta,w)
\in ^{}V^{2}\times V$
such that
\begin{equation}
\int_{\Omega}  D^* \kappa(\underline
\theta) :  \kappa(\underline{\hat\theta})\;d\Omega+ kG\int_\Omega (\nabla w - \underline \theta) \cdot (\nabla
\hat w - \underline{\hat\theta})\;d\Omega = \lambda \int_\Omega  \sigma\nabla w \cdot \nabla \hat w\;d\Omega
\quad \forall (\underline{\hat\theta},\hat w) \in V\times V.
\label{eqn:MomentWeakForm}
\end{equation}

\smallskip\noindent
We approximate the solutions of \eqref{eqn:MomentWeakForm} using again a piecewise bilinear finite elements on a quadrilateral mesh $\mathcal Q_h$,  and such that
$w$ and $\underline \theta$ are interpolated with the same shape
functions $\{\phi_i(x,y)\}^{\rm{nnod}}_{i=1}$. The matrix form of
\eqref{eqn:MomentWeakForm} is 
\begin{equation}
\textbf K^{\rm{B}} \underline d^{\rm{B}} = \lambda \textbf G\underline d^{\rm{B}},
\label{eqn:matrixeig}
\end{equation}
where $\textbf  K^{\rm{B}}\in \mathbb{R}^{M \times M}$ is the global
stiffness matrix  (LHS of \eqref{eqn:MomentWeakForm}) whilst $\textbf
G \in \mathbb R^{M \times M}$ is the geometric stiffness matrix (RHS
of \eqref{eqn:MomentWeakForm}). Further details of the exact finite element formulation, for a similar eigenvalue problem are provided in \cite{Dod16}.

\subsection{Results I: Comparison between MC, MLMC and MLMC-SR}\label{sec:results}

\begin{figure}
\centering
\includegraphics[width = 0.18\linewidth]{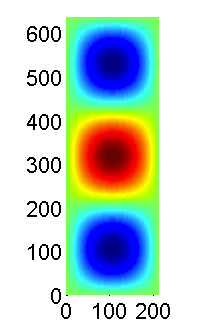}{\includegraphics[width = 0.38\linewidth]{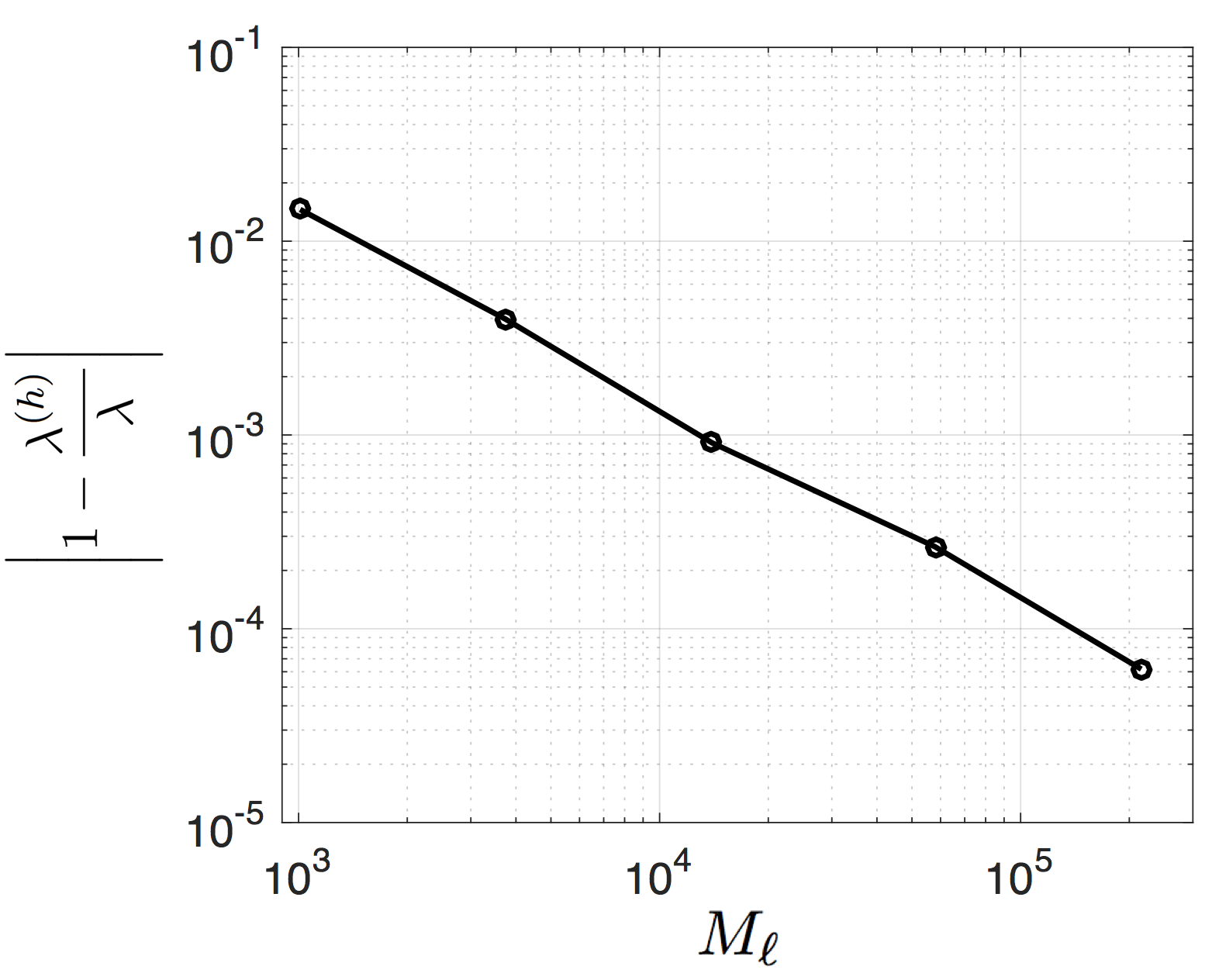}}\quad{\includegraphics[width = 0.38\linewidth]{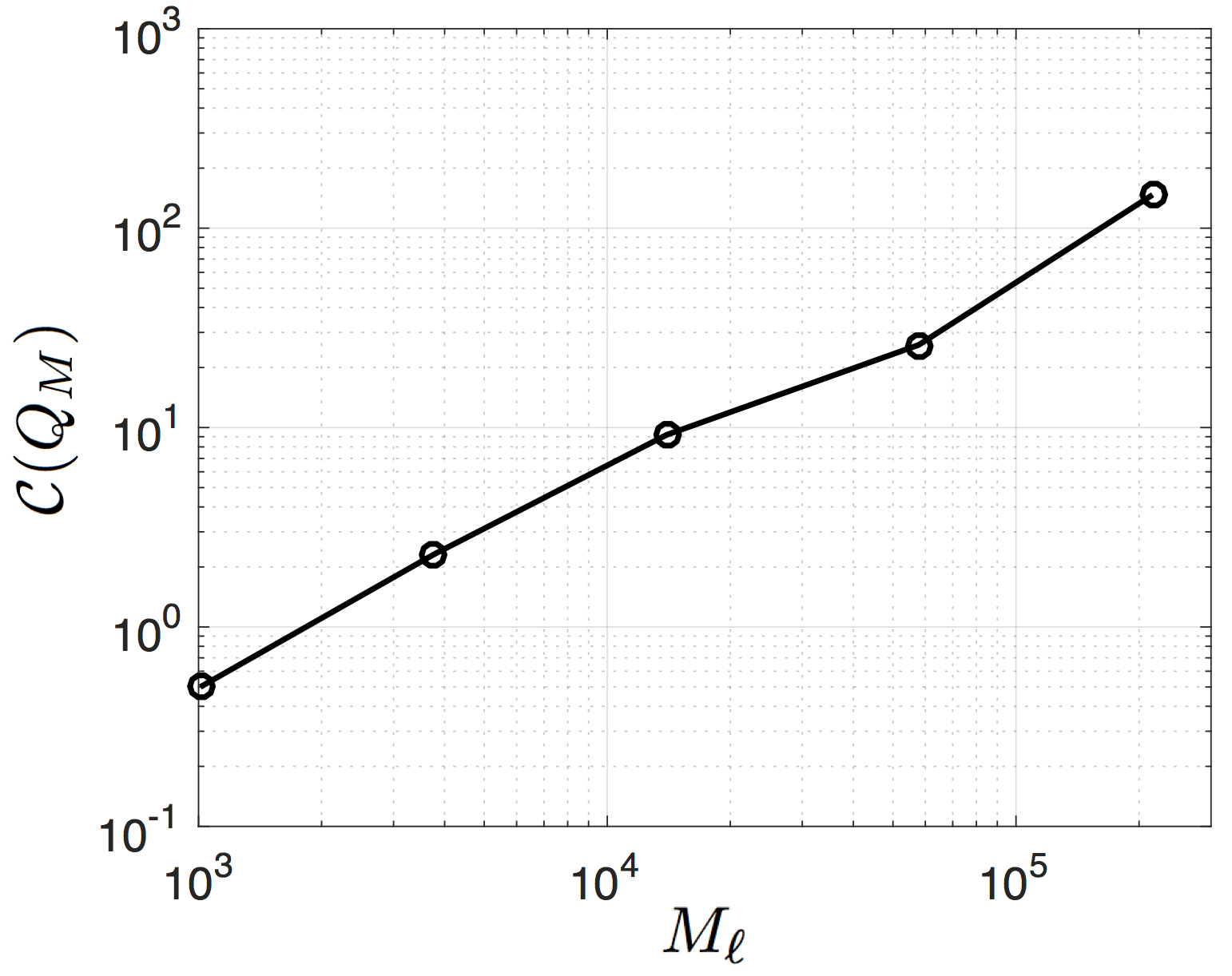}}
\caption{(Left) Plot of the critical buckling mode of the pristine panel corresponding to the critical buckling load of $278.59$kN. (Middle) Log-Log plot of the relative FE error in the buckling load $|1 - \lambda^{(h)}/\lambda|$ against $M$, which shows that the error converges with order $\alpha \approx 1$ (Right) Log-Log plot of Cost (CPU-time) against $M$ showing $\mathcal{C}(Q_M) \simeq M^{1.17}$ (i.e $\gamma = 1.17$).}
\label{fig:pristine}
\end{figure}

Before comparing MC, MLMC and MLMC-SR for the second test problem, we first consider the convergence rates for the FE approximation of the critical buckling load $\lambda$, as well as the associated computational cost $(\rm{Cost})$ under uniform mesh refinement. Figure \ref{fig:pristine} (middle) shows the convergence of the relative error in $\lambda$. We see that for the pristine case, the buckling load converges to a value of $278.59$kN (the mode is shown), at a rate $\alpha \approx 1$ with respect to the number of degrees of freedom $M$, i.e. 
$$
\left| 1- \frac{\lambda^{(h)}}{\lambda}\right| \leq CM^{-1},
$$
for some constant $C>0$, independent of the number of degrees of freedom $M$. This agrees with the theoretically predicted convergence rate for buckling modes for this element. We approximate the value for $\gamma$, the rate at which the Cost (in CPU-time) scales with $M$, as shown in Fig. \ref{fig:pristine} (right). The gradient of the line shows that $$\mathcal{C}(Q_M) \leq C M^{1.17},$$ i.e. a value of $\gamma \approx 1.17$. The CPU-time is made up of matrix assembly for \eqref{eqn:matrixeig} and the calculation of the smallest eigenvalue of \eqref{eqn:matrixeig}. For the size of problems considered here ($\ell \leq 8$), the CPU-time is dominated by the matrix assembly, which scales linearly with $M$. For larger problem sizes ($M \simeq 1e^6$), the eigenvalue solve will dominate the CPU-time and $\gamma$ will increase as the limit of eigenvalue solvers (\texttt{ARPACK}) for 2D problems is reached.

\smallskip
\noindent
In this test we estimate the $\mathbb P(\lambda < \lambda^\star = 272.47\mbox{kN})$ this corresponds to estimating the mean value of the quantity of interest
\begin{equation}
Q(\lambda) = 
\begin{cases}
1 \quad \mbox{if} \quad \lambda < \lambda^\star \\
0 \quad \mbox{if} \quad \lambda \geq \lambda^\star
\end{cases}
\end{equation} 
In the results the initial refinement level $\ell = 0$ is created by five uniform refinements mesh containing a single element; this is to ensure that some failures are observed on the coarsest level, and is necessary due to the one-sided convergence of the buckling load (i.e. buckling load only reduces with mesh refinement).

\begin{figure}\label{fig:MLMC_SR}
\centering
\includegraphics[width = 0.9\linewidth]{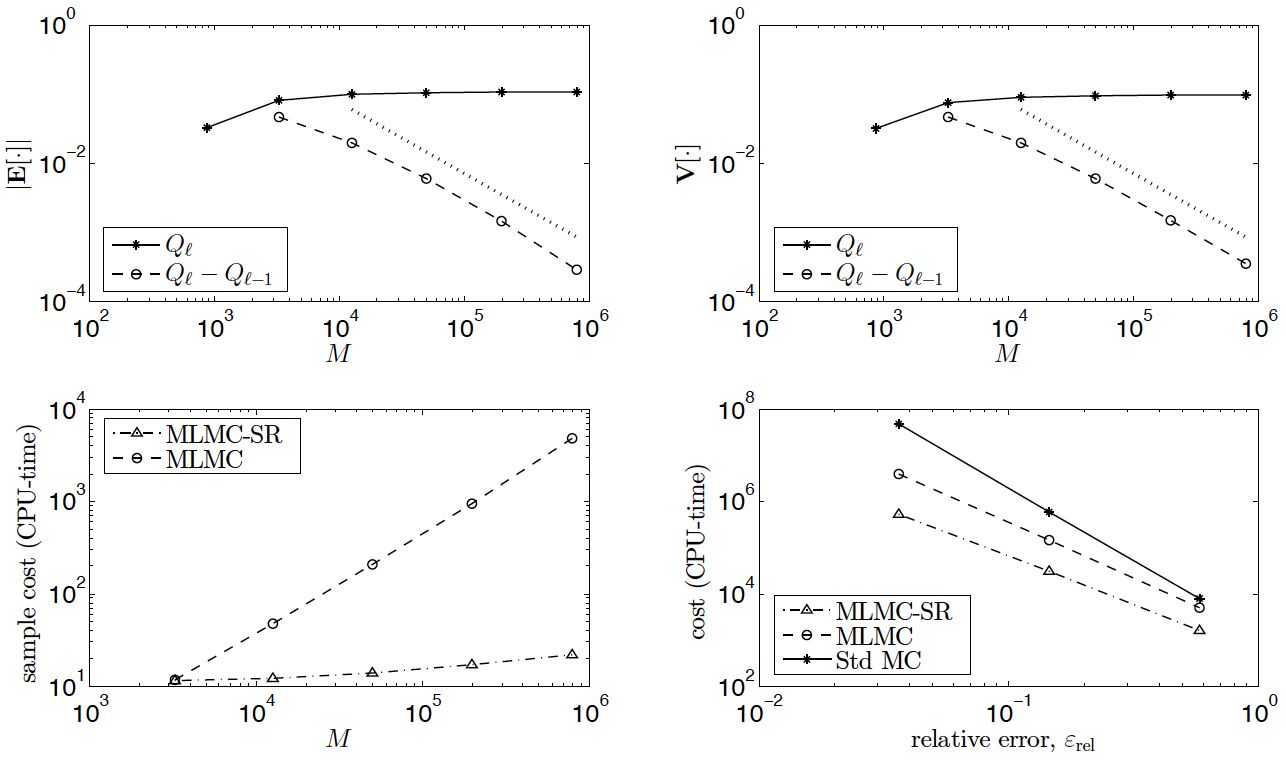}
\caption{(Top-Left) Expected values of $Q_\ell$ and $Y_\ell$ against degrees
of freedom $M$. The gradient of the dotted line is $1.03$. (Top-Right) Variance of $Q_\ell$ and $Y_\ell$ against degrees of freedom $M$. The gradient of the dotted line is $1.03$. (Bottom-Left) Comparison of expected cost per sample (CPU-time)
for MLMC-SR and MLMC. The gradient of the dash-dotted line (MLMC-SR) is $0.12$, whilst the
gradient of the dashed line (MLMC) is $1.16$. (Bottom-Right) Comparison of cost (CPU-time) for
MLMC-SR, MLMC and standard MC, against relative error tolerance $e_{rel}$. The gradient of the
dash-dotted line (MLMC-SR) is $2.03$, the gradient of the dashed line (MLMC) is $2.28$, whilst
the gradient of the solid line (MC) is $3.14$.}
\end{figure}

\smallskip
The MLMC-SR simulation has been carried out for an error tolerance of
$e = 4.00 \times 10^{-3}$,
corresponding to an approximate relative error $3.6\%$. Again we split the error equally between bias and
sampling error ($\theta=1/2$). Figure \ref{fig:MLMC_SR} (top-left and top-right) shows the behaviour of the expected
value and variance of $Y_\ell$, with respect to degrees of freedom $M_\ell$. From this we estimate that
\begin{equation}
\mathbb E[Y_\ell] \simeq M_\ell^{-1.03} \quad \mbox{and} \quad \mathbb V[Y_\ell] \simeq M_\ell^{-1.03}
\end{equation}
and hence $\alpha \approx \beta \approx 1.03$. This is in agreement
with the theoretically predicted convergence rates.

\smallskip\noindent The bottom two plots in Figure \ref{fig:MLMC_SR} show the computational cost of the MLMC-SR simulation, in comparison with the standard MC and MLMC simulations. The lower-left plot compares the expected cost per sample for the MLMC-SR and MLMC simulations, with respect to degrees of freedom. From \eqref{eqn:costpersample}, taking parameter values $\alpha = 1.03$ and $\gamma=1.17$ (as approximated above), the predicted growth rate for MLMC-SR is $M_\ell^{\gamma - \alpha} = M_\ell^{0.14}$. This is verified by the numerical results, for which we observe a growth rate of $M_\ell^{0.12}$. For MLMC we observe the significantly greater cost growth of approximately $M_{\ell}^{1.16}$.

\smallskip\noindent
The lower-right plot shows the computational cost of the MLMC-SR simulation versus that of MLMC and standard MC, for a range of relative error tolerances. In this regime, we have $\alpha < \gamma < 2\alpha$, for which we predict the cost MLMC-SR simulation to grow proportionally
to $\tau^{-2}$. This agrees with an observed rate $\tau^{-2.03}$. The costs of the MLMC and standard MC
simulations are predicted to grow like $\tau^{-2.13}$ and $\tau^{-3.13}$, respectively, and again there is good
agreement with the numerical approximations $\tau^{-2.28}$, and $\tau^{-3.14}$.
Table 3.5 lists the optimal numbers N` of samples required by the MLMC-SR
and MLMC simulations for these error tolerances. The total computational costs of the simulations
are included, along with the corresponding costs for the standard MC simulation. For the smallest
error tolerance considered, $3.6\%$, the MLMC-SR simulation reduces the computational cost by a factor of $7.43$ compared to the MLMC simulation, and $90.32$ as compared to the
standard MC simulation. 

\smallskip\noindent
The distribution of work across refinement levels for the MLMC-SR simulation (in the case $\tau_{rel.} = 3.6\%$) is presented in Table \ref{tab:MLMCSR_work}. As expected, we observe that very few realisations are solved on their highest refinement levels, with most of the computational effort restricted to the coarser levels.

\begin{table}[]\label{tab:MLMCSR_work}
\begin{tabular}{|c|c|c|c|c|c|c|c|c|c|}
\hline
\multirow{2}{*}{$\tau_{rel}$} & \multirow{2}{*}{Method} & \multicolumn{6}{c|}{$N_\ell$}                     & \multirow{2}{*}{Cost(hrs)} & \multirow{2}{*}{Saving} \\ \cline{3-8}
                                     &                         & 0      & 1      & 2      & 3     & 4     & 5      &                             &                                   \\ \hline
\multirow{3}{*}{3.6\%}               & {\small MLMC-SR}                 & 26,883 & 16,489 & 10,348 & 5,345 & 2,407 & 1,029  & 147                         & -                                 \\
                                     & {\small MLMC }                  & 73,226 & 44,365 & 14,382 & 3,754 & 879   & 189    & 1,097                       & 7.43                              \\
                                     & {\small MC }                    & -      & -      & -      & -     & -     & 11,897 & 13,341              & 90.32                             \\ \hline
\multirow{3}{*}{15\%}                & {\small MLMC-SR }                & 1,610  & 988    & 620    & 320   & 145   & -      & 8.50                        & -                                 \\
                                     & {\small MLMC }                   & 3,523  & 2,135  & 692    & 181   & 43    & -      & 41                          & 4.81                              \\
                                     & {\small MC }                     & -      & -      & -      & -     & 744   & -      & 162                         & 19.24                             \\ \hline
\end{tabular}
\caption{Comparison of optimal number of samples $N_\ell$ and computational cost for MLMC-SR, MLMC, and standard MC simulations for Test Problem II.}
\end{table}

\subsection{Estimation of rare events}

\smallskip\noindent
We now push the MLMC-SR methodology to estimate a much smaller failure probability of approximately $\sim 1 / 150 = 0.00\dot{6}$. This will push the computational demand well beyond the reach of standard Monte Carlo, and demonstrate the potential computational benefits of adopting the multilevel strategies for the estimation of rare events. In these numerical experiments we use the same setup as described above and estimate
$$
\mathbb P (\lambda \leq \lambda^\star = 268\mbox{kN}).
$$
However, the standard multilevel approach gives rise to somewhat of a paradox. If we wish to estimate a rare event $\mathbb E[Q]$, in the multilevel framework we adopt the multilevel telescoping sum \eqref{eq:MLMC_expectation}. For higher levels this leaves us trying to estimate $\mathbb E[Y_\ell]$. In our stochastic eigenvalue examples $Y_\ell = 1$ only if failure occurs on level $\ell$ but not on level $\ell-1$. Of course, this conditional probability, is a much rarer event than failure occurring on level $\ell$. Paradoxically for the higher levels, we would require many more simulations to see just one case where the two adjacent levels differ. However, we note that with the use of selective refinement in most cases samples are pre-screened by coarser/cheaper model solves; so only very rarely, when there is a discrepancy at higher levels, do we require expensive solves. Mostly, for a rare event, initial coarse solves are sufficiently far away from $\lambda^\star$ to guarantee that the particular sample does not fail on any level according to \eqref{eqn:testFailure}. As a result, in the calculations we present below, the average computational cost of a sample on level $2$ is little different to that on level 5, $0.239$ secs compared with $0.244$ secs respectively.

\smallskip
This highlights that in such cases significant computational gains can be achieved by using MLMC-SR  with a simple two level multilevel estimate, i.e.
$$
\mathbb E[Q] \approx \hat Q_0 +  \hat Y_{L,0}, \quad \mbox{where} \quad Y_{i,j} = Q_i - Q_j.
$$
Importantly we note here, that we still use all the levels of refinement to calculate $Q_L$ in the selective refinement procedure. We then simply only use the coarse and fine levels in the MLMC estimate. Because of the plateau in cost for MLMC-SR for rare events, it is less efficient to use all levels in the estimator, yet we can still exploit some variance reduction with a two level method. This fact, highlights that in the results to follow most of the computational gains come from the selective refinement strategy.

\smallskip \noindent
Firstly, for our tightest tolerance we estimate that 
$$
\mathbb P (\lambda \leq \lambda^\star = 268\mbox{kN}) = 0.00645.
$$
Table \ref{tab:rare} summarises the computational savings of MLMC-SR over standard MC, over a range of tolerances. In each case bias and sampling error are balanced (i.e $\theta = 1/2$). We note that because of the scale of these calculations results were computed on \texttt{Isca}, Exeter's supercomputer $\sim 400$ nodes each with 8-core Intel Xeon E5-2650v2 Ivybridge processors each running at 2.6 GHz and giving a total of over 6000 available cores. Our calculations were trivially distributed over 1,024 processors. We see that over the range of tolerances MLMC-SR demonstrates huge computational savings. In particular for our finest tolerance calculation, we see an estimated saving of a factor of $1173$. Most importantly, the scale of computation required from standard MC would require $218$ days of computation on a large computing resource, the MLMC-SR reduces this to just a few hours. Importantly from an engineering perspective, this scale of savings opens the opportunity to new studies of rare events.

\begin{table}[]
\label{tab:rare}
\begin{tabular}{|c|c|c|c|c|c|c|c|c|c|}
\hline
\multirow{2}{*}{$\tau_{rel}$} & \multirow{2}{*}{Method}   & \multirow{2}{*}{Term} & \multicolumn{5}{c|}{Solves on level}                                                                                                                             & \multirow{2}{*}{Cost} & \multirow{2}{*}{Saving} \\ \cline{4-8}
                            &                           &                                     & 0                                                & 1                                                & 2                & 3                  & 4                  &                                           &                                        \\ \hline
\multirow{3}{*}{4.3\%}      & \multirow{2}{*}{MLMC-SR}  & $\hat{Q}_0$                    & 3.65e5                                          & -                                                & -                & -                  & -                  & \multirow{2}{*}{$35.7$secs}               & \multirow{2}{*}{69}                    \\
                            &                           & $\hat{Y}_{2,0}$              & 2.54e5                                          & 2.53e5                                        & 348              & -                  & -                  &                                           &                                        \\ \cline{2-10} 
                            & MC                        & $\hat Q_2$                    & -                                                & -                                                & 3.18e5        & -                  & -                  & $41$mins                                 & -                                      \\ \hline
\multirow{3}{*}{1.4\%}      & \multirow{2}{*}{MLMC-SR} & $\hat{Q}_0$                    & $3.32e6$ & -                                                & -                & -                  & -                  & \multirow{2}{*}{$5.46$mins}              & \multirow{2}{*}{124}                   \\
                            &                           & $\hat{Y}_{3,0}$              & $2.41e6$ & $2.40 e6$ & 4,268            & 965                & -                  &                                           &                                        \\ \cline{2-10} 
                            & MC                        & $\hat Q_3$                    & -                                                & -                                                & -                & $2.93e6$ & -                  & $11.24$hrs                               & -                                      \\ \hline
\multirow{3}{*}{0.2\%}      & \multirow{2}{*}{MLMC-SR} & $\hat{Q}_0$                    & $1.63 e8$ & -                                                & -                & -                  & -                  & \multirow{2}{*}{$4.4$hrs}              & \multirow{2}{*}{1173}                  \\
                            &                           & $\hat{Y}_{4,1}$              & $1.19 e8$                               & $1.18 e8$                               & $2.19 e5$ & $5.61e4$   & $8,348$            &                                           &                                        \\ \cline{2-10} 
                            & MC$^\star$                       & $\hat Q_4$                    & -                                                & -                                                & -                & -                  & $1.44e8$ & $218$days                                & -     
                            \\ \hline                                
\end{tabular}
\caption{Demonstrates the relative saving of MLMC-SR over MC for the estimation of a rare event for a range of tolerances, alongside the distribution of work on each level for MLMC-SR. The cost is that of a simulation distributed over $1024$ processors ($\star$ indicates a calculation which can only be estimated due to scale of computation). }
\end{table}

\section{Conclusions}

In this paper we have successfully demonstrated the applicability of MLMC simulation on two typical aerospace model problems. From our numerical results, the advantages of MLMC simulation over standard MC simulation are apparent, with huge savings in computational cost being observed. We see also that MLMC simulation is not limited to easy problems, and in fact the gains are more pronounced in cases where the discretisation error is large.  We have further demonstrated the versatility of MLMC simulation, showing that the method is not restricted to problems in which the quantity of interest is a smooth functional of the solution vector, but can readily be applied and extended to calculated failure probabilities with significant computational speed-ups.

\smallskip
From an engineering viewpoint, whilst the model problems are chosen to represent the typical gains achieved by the MLMC methodology, in addition, we learn
something about the engineering implications of uncertainty in each case. In the buckling test problem, perhaps unsurprisingly, the numerical results show that random
variations in ply angles increase the risk of buckling failure significantly. With ply angles of the order typically observed in an Automated Fibre Placement (AFP) machine ($\pm 5^\circ$) significant variability is observed in buckling performance. As for our numerical results into the effects of random fibre waviness on the compressive strength of composites, high fidelity stochastic simulations show a remarkable agreement with Budiansky's classical kinking model \cite{Bud83} if the misalignment angle is taken to be the standard deviation of the misalignment random field.

\smallskip
Current and future research is exploring the use of sample-dependent adaptive grids, to exploit the computational gains offered by adaptive finite elements \cite{Det18}; as well as integrating the Multilevel Framework with experimental data in a Bayesian setting to quantify and reduce modelling uncertainties as proposed by theoretical methodology introduced in Dodwell et al. (2015) \cite{Dod15}.


\section*{Acknowledgements}

This work falls within EPSRC Project EP/K031368/1
  ``Multiscale Modelling of Aerospace Composites''. Dodwell was supported by The Alan Turing Institute under the EPSRC Grant EP/N510129/1 and Butler holds a Royal Academy of Engineering-GKN Aerospace Research Chair in Composites.

\end{document}